\documentclass[a4paper,12pt]{amsart}
\usepackage{ifthen}
\usepackage{mathrsfs}
\numberwithin{equation}{section}
\nonstopmode
\newtheorem{thm}{Theorem}[section]
\newtheorem{cor}[thm]{Corollary}
\newtheorem{lem}[thm]{Lemma}
\newtheorem{prop}[thm]{Proposition}

\newtheorem{prob}[thm]{Problem}

\theoremstyle{definition}

\newenvironment{rem}{%
\bigskip
\noindent
\textsl{{\sl Remark. }}}{\bigskip}

\newenvironment{pf}[1][]{%
 \vskip 3mm
 \noindent
 \ifthenelse{\equal{#1}{}}%
  {{\slshape Proof. }}%
  {{\slshape #1.} }%
 }%
{\qed\bigskip}

\newcounter{alphabet}
\newcounter{tmp}
\newenvironment{Thm}[1][]{\refstepcounter{alphabet}%
\bigskip%
\noindent%
{\bf Theorem \Alph{alphabet}}%
\ifthenelse{\equal{#1}{}}{}{ (#1)}%
{\bf .}
\itshape}{\vskip 8pt}

\newcommand{\A}{{\mathcal A}}

\newcommand{\C}{{\mathbb C}}
\newcommand{\D}{{\mathbb D}}

\renewcommand{\Re}{{\operatorname{Re}\,}}

\newcommand{\inv}{^{-1}}

\newcommand{\arctanh}{{\operatorname{arctanh}\,}}

\newcommand{\aand}{{\quad\text{and}\quad}}

\newcounter{minutes}
\divide\time by 60
\newcounter{hours}
\multiply\time by 60
\addtocounter{minutes}{-\time}
\begin{document}
\bibliographystyle{amsplain}
\title{
Univalence criteria and analogs of the John constant
}

\author[Y.~C.~Kim]{Yong Chan Kim}
\address{Department of Mathematics Education, Yeungnam University, 214-1 Daedong
Gyongsan 712-749, Korea}
\email{kimyc@ynu.ac.kr}
\author[T. Sugawa]{Toshiyuki Sugawa}
\address{Graduate School of Information Sciences,
Tohoku University, Aoba-ku, Sendai 980-8579, Japan}
\email{sugawa@math.is.tohoku.ac.jp}
\begin{abstract}
Let $p(z)=zf'(z)/f(z)$ for a function $f(z)$ analytic on the unit disk
$|z|<1$ in the complex plane and normalized by $f(0)=0, f'(0)=1.$
We will provide lower and upper bounds for the best constants
$\delta_0$ and $\delta_1$ such that the conditions
$e^{-\delta_0/2}<|p(z)|<e^{\delta_0/2}$ and $|p(w)/p(z)|<e^{\delta_1}$
for $|z|,|w|<1$ respectively imply univalence of $f$ on the unit disk.
\end{abstract}
\thanks{
The first author was supported by Yeungnam University (2012).
The second author was supported in part by JSPS Grant-in-Aid for 
Scientific Research (B) 22340025.
}
\maketitle

\section{Introduction}
For a non-constant analytic function $f$ on the unit disk $\D=\{z\in\C:|z|<1\},$
set
$$
M(f)=\sup_{z\in\D}|f'(z)|
\aand
m(f)=\inf_{z\in\D}|f'(z)|.
$$
Note that $M(f)$ is a positive number (possibly $+\infty$)
whereas $m(f)$ is a finite nonnegative number.
F.~John \cite{John69} proved the following result.

\begin{Thm}[John (1969)]
There exists a number
$\gamma\in[\pi/2,~\log(97+56\sqrt3)]$ with the following
property:
if a non-constant analytic function $f$ on $\D$ satisfies the condition
$M(f)\le e^\gamma m(f),$ then $f$ is univalent on $\D.$
\end{Thm}

We remark that $\log(97+56\sqrt3)=5.2678\dots.$
The largest possible number $\gamma$ with the property in the theorem
is called the (logarithmic) John constant and will be denoted by $\gamma_1.$
(In the literature, the John constant refers to $e^{\gamma_1}.$ 
We adopt, however, the logarithmic one for our convenience in this note.)
Yamashita \cite{Yam78} improved John's result by showing that $\gamma_1\le \pi.$
Gevirtz \cite{Gev81, Gev89} further proved that $\gamma_1\le \lambda \pi$
and conjectured that $\gamma_1=\lambda\pi,$
where $\lambda=0.6278\dots$ is the number determined by
a transcendental equation.

We could consider a similar problem for $zf'(z)/f(z)$ instead of $f'(z)$
for an analytic function $f$ on $\D$ with $f(0)=0, f'(0)\ne0.$
Let
$$
L(f)=\sup_{z\in\D}\left|\frac{zf'(z)}{f(z)}\right|
\aand
l(f)=\inf_{z\in\D}\left|\frac{zf'(z)}{f(z)}\right|
$$
for such a function $f.$
Here, the value of $zf'(z)/f(z)$ at $z=0$ will be understood as
$\lim_{z\to0}zf'(z)/f(z)=1$ as usual.
Note that $0\le l(f)\le1\le L(f)\le+\infty.$
It is easy to see that $l(f)=1$ (or $L(f)=1$)
precisely when $f(z)=az$ for a nonzero constant $a.$
Since $zf'(z)/f(z)$ is unchanged under the dilation $f\mapsto af$ for a nonzero
constant $a,$ we can restrict our attention to analytic functions $f(z)$
on $\D$ normalized by $f(0)=0, f'(0)=1.$
The class of those normalized analytic functions on $\D$
will be denoted by $\A$ in the sequel.
Thus the problem can be formulated as follows.

\begin{prob}\label{prob}
Find a number $\delta>0$ with the following property:
If a function $f\in\A$ satisfies the condition
$L(f)\le e^\delta l(f)$ then $f$ is univalent on $\D.$
\end{prob}

Since the value $1$ plays a special role in the study of $zf'(z)/f(z),$
the following problem is also natural to consider.

\begin{prob}\label{prob0}
Find a number $\delta>0$ with the following property:
If a function $f\in\A$ satisfies the condition
$e^{-\delta/2}<|zf'(z)/f(z)|<e^{\delta/2}$ on $\D,$ then
$f$ is univalent on $\D.$
\end{prob}

Let $\delta_1$ and $\delta_0$ be the largest possible numbers $\delta$
in Problems \ref{prob} and \ref{prob0}, respectively (if they exsist).

The authors proved in \cite{KS12PT} that 
$\pi/6=0.523\dots\le\delta_0\le\pi=3.14\dots.$
Obviously, $\delta_1\le\delta_0\le 2\delta_1.$
Therefore, we already have the estimates $\pi/12=0.261\dots\le\delta_1\le\pi.$

The purpose of the present note is to improve the estimates.

\begin{thm}\label{thm:main1}
$\dfrac{\pi}3=1.04719\dots<\delta_0<\dfrac{5\pi}{7}=2.24399\dots.$
\end{thm}

\begin{thm}\label{thm:main2}
$\dfrac{7\pi}{25}=0.87964\dots<\delta_1<\dfrac{5\pi}{7}=2.24399\dots.$
\end{thm}

We remark that the above results are not optimal.
Indeed, more elaborative numerical computations would yield
slightly better bounds as will be suggested at the end of Section 2.

In order to obtain a lower bound, we need a univalence criterion
due to Becker \cite{Becker72}
with numerical computations as we will explain in Section 2.
On the other hand, to give an upper bound, we should construct a non-univalent
function satisfying the condition in Problems \ref{prob} or \ref{prob0}.
The function $F_a\in\A$ determined by
the differential equation
\begin{equation}\label{eq:Fa}
\frac{zF_a'(z)}{F_a(z)}
=\left(\frac{1-iz}{1+iz}\right)^{ai}
\end{equation}
is a candidate for an extremal one, where $a$ is a positive constant
and $i$ is the imaginary unit $\sqrt{-1}.$
As will be seen later, $L(F_a)/l(F_a)=e^{\pi a}.$
We will give a detailed account on this function and provide the upper
bound in the above theorems in Section 3.
The proof is involved with matrices of large order.
Therefore, we made use of Mathematica 8.0 to carry out symbolic computations.

The most interesting problem is to determine the values of $\delta_0$
and $\delta_1.$
However, this seems to the authors very hard.
We end the section with a couple of open questions, which may be
easier to solve.
Let $a^*$ be the supremum of the numbers $a$ such that $F_a$ is univalent
on $\D.$
Likewise let $a_*$ be the infimum of the numbers of $a$ such that
$F_a$ is not univalent on $\D.$
Obviously, $\delta_0\le\pi a_*\le\pi a^*.$
In the proof of the above theorems, we indeed show that $a_*<5/7.$

\begin{enumerate}
\item[(1)] Is it true that $a_*=a^*?$
\item[(2)] Is it true that $\delta_0=\pi a_*?$
\item[(3)] Is it true that $\delta_0=\delta_1?$
\end{enumerate}

\medskip
\noindent
{\bf Acknowledgement.}
The authors would like to thank the referee for careful reading of the
manuscript and for suggestions which helped us to improve the exposition.

\section{Obtaining lower bounds: univalence criteria}

We recall basic hyperbolic geometry of the unit disk $\D.$
The hyperbolic distance between two points $z_1,z_2$ in $\D$ is defined by
$$
d(z_1,z_2)=\inf_\gamma\int_\gamma\frac{|dz|}{1-|z|^2},
$$
where the infimum is taken over all rectifiable paths $\gamma$
joining $z_1$ and $z_2$ in $\D.$
The Schwarz-Pick lemma asserts that
\begin{equation}\label{eq:sp}
\frac{|\omega'(z)|}{1-|\omega(z)|^2}\le \frac1{1-|z|^2},\quad z\in\D,
\end{equation}
for any analytic map $\omega:\D\to\D.$
In particular, for an analytic automorphism $T$ of $\D,$
we have $|T'(z)|/(1-|T(z)|^2)=1/(1-|z|^2)$ and therefore,
$d(T(z_1), T(z_2))=d(z_1,z_2)$ for $z_1,z_2\in\D.$
It is well known that the above infimum is attained by the circular arc
(possibly a line segment)
joining $z_1$ and $z_2$ whose whole circle is perpendicular to the unit circle.
By using these facts, one can compute the hyperbolic distance:
$d(z_1,z_2)=\arctanh|(z_1-z_2)/(1-\bar z_1 z_2)|.$
Here, we recall that $\arctanh r=\frac12\log\frac{1+r}{1-r}.$

The following is a useful univalence criterion due to Becker \cite{Becker72}.

\begin{lem}\label{lem:becker}
Let $f$ be a non-constant analytic function on $\D.$
If
$$
(1-|z|^2)\left|\frac{zf''(z)}{f'(z)}\right|\le 1,\quad z\in\D,
$$
then $f$ is univalent on $\D.$
\end{lem}

Sometimes, it is more convenient to consider the pre-Schwarzian norm
$$
\|f\|=\sup_{z\in\D}(1-|z|^2)\left|\frac{f''(z)}{f'(z)}\right|
$$
because it has several nice properties (see \cite{KPS02} for example).
By Becker's theorem above, we see that the condition $\|f\|\le 1$
implies univalence of $f$ on $\D.$
We used this norm to deduce the estimate $\pi/6\le\delta_0.$
In this note, however, we do use the original form (Lemma \ref{lem:becker})
of Becker's theorem to improve the estimate.

For a non-negative number $c,$ we consider the quantity
$$
\Phi(c)
=\sup_{0<r<1}\Big\{r+c(1-r^2)\arctanh r\Big\}
=c\sup_{0<r<1}\Big\{c\inv r+(1-r^2)\arctanh r\Big\}.
$$
It is easy to see that $\Phi(c)$ is non-decreasing in $c$ and that
$c\inv\Phi(c)$ is non-increasing in $c.$
In terms of this function, we will prove the following technical lemma
which yields lower bounds for $\delta_0$ and $\delta_1$ as corollaries.

\begin{lem}\label{lem:key}
Let $f\in\A.$
If $L(f)/l(f)<+\infty$ and if the inequality
$$
\frac{2}{\pi}\Phi(L(f))\log\frac{L(f)}{l(f)}\le 1
$$
holds, then $f$ is univalent on $\D.$
\end{lem}

The lemma immediately yields the following results.

\begin{cor}\label{cor:main}
Let $\delta>0$ be given.
If
\begin{equation}\label{eq:0}
\frac{2\delta}\pi \Phi(e^{\delta/2})\le1,
\end{equation}
then $\delta\le\delta_0.$
If
\begin{equation}\label{eq:1}
\frac{2\delta}\pi \Phi(e^{\delta})\le1,
\end{equation}
then $\delta\le\delta_1.$
\end{cor}

To show the corollary, we first assume \eqref{eq:0} and
consider a function $f\in\A$ satisfying
$e^{-\delta/2}<|zf'(z)/f(z)|<e^{\delta/2}.$
Then $L(f)\le e^{\delta/2}$ and $\log L(f)/l(f)\le\delta$ so that
$$
\frac{2}{\pi}\Phi(L(f))\log\frac{L(f)}{l(f)}
\le\frac{2\delta}\pi \Phi(e^{\delta/2})\le1.
$$
We now apply Lemma \ref{lem:key} to conclude univalence of $f.$
Secondly, we assume \eqref{eq:1} and consider a function $f\in\A$
satisfying $L(f)\le e^\delta l(f).$
Then $L(f)\le e^\delta$ and the conclusion follows similarly.

Let us prepare for the proof of Lemma \ref{lem:key}.
We note that the function $\arctan z=\frac1{2i}\log\frac{1+iz}{1-iz}$ 
maps  the unit disk $\D$ conformally onto
the vertical parallel strip $|\Re w|<\pi/4.$ 
Therefore, for a constant $a>0,$ the function
\begin{equation}\label{eq:Q}
Q_a(z)=\exp(2a\arctan z)
%=\exp(ai\log\tfrac{1-iz}{1+iz})
=\left(\frac{1-iz}{1+iz}\right)^{ai}
\end{equation}
is the universal covering projection of $\D$
onto the annulus $e^{-\pi a/2}<|w|<e^{\pi a/2}.$
We note that the function $Q_a$ satisfies $Q_a(0)=1$ and
$$
\frac{Q_a'(z)}{Q_a(z)}=\frac{2a}{1+z^2}.
$$

\begin{pf}[Proof of Lemma \ref{lem:key}]
Let $p(z)=zf'(z)/f(z)$ for a function $f\in\A.$
If $p$ is a constant, then $f$ is clearly univalent.
We can thus assume that $p$ is not a constant so that $l(f)<1<L(f).$
Let $\delta=\log L(f)/l(f)<+\infty$ and $m=\sqrt{L(f)l(f)}.$
We consider the universal covering map $Q=mQ_a$
of $\D$ onto the annulus $W=\{w: l(f)<|w|<L(f)\}
=\{w: me^{-\delta/2}<|w|<me^{\delta/2}\},$
where $Q_a$ is given in \eqref{eq:Q} with $a=\delta/\pi.$
Note that $p(\D)\subset W$ by assumption.
Since the real interval $(-1,1)$ is mapped onto $(l(f), L(f))$
by $Q,$ we can choose an $\alpha\in(-1,1)$ so that $Q(\alpha)=1.$
Then, $P=Q\circ T$ is a universal covering map of $\D$ onto $W$
with $P(0)=1,$ where $T(z)=(z+\alpha)/(1+\alpha z).$ 
Since $P:\D\to W$ is a covering map, we can take a lift $\omega$ of
$p$ with respect to $P$ so that $\omega(0)=0$ and $p=P\circ\omega.$
We write $w=\omega(z).$
Note here that the Schwarz lemma implies $|w|\le |z|.$
We now have
$$
\frac{zf''(z)}{f'(z)}
=\frac{zp'(z)}{p(z)}+p(z)-1
=\frac{z\omega'(z)P'(w)}{P(w)}+P(w)-1.
$$
Set $\tau=T(w)\in\D.$
Since $T$ is a hyperbolic isometry of $\D,$
one has the relation $(1-|w|^2)|T'(w)|=1-|\tau|^2.$
Therefore, by using \eqref{eq:sp}, we have
\begin{align*}
(1-|z|^2)\left|\frac{\omega'(z) P'(w)}{P(w)}\right|
&\le (1-|w|^2)\left|\frac{Q'(\tau)T'(w)}{Q(\tau)}\right| \\
&=(1-|\tau|^2)\left|\frac{Q'(\tau)}{Q(\tau)}\right| \\
&=\frac{2a(1-|\tau|^2)}{|1+\tau^2|} \\
&\le 2a.
\end{align*}
Let $\gamma$ be the image of the line segment $(0,w)$ under
the M\"obius mapping $T.$
Then,
$$
P(w)-1
=\int_0^w P'(t)dt 
=\int_0^wQ'(T(t))T'(t)dt 
=\int_\gamma Q'(u)du 
=\int_\gamma \frac{2aQ(u)}{1+u^2}du.
$$
Since $|Q(u)|\le L(f),$ we obtain
\begin{align*}
|P(w)-1|
&\le 2aL(f)\int_\gamma\frac{|du|}{1-|u|^2} 
= 2aL(f)\int_0^w\frac{|du|}{1-|u|^2} \\
%&=2aL(f) d(T(0), T(w)) \\
&=2aL(f) d(0, w) 
\le 2aL(f)\arctanh |z|.
\end{align*}
Therefore, 
\begin{equation}\label{eq:ps}
(1-|z|^2)\left|\frac{zf''(z)}{f'(z)}\right|
\le 2a|z|+2aL(f)(1-|z|^2)\arctanh |z|.
\end{equation}
Hence,
$$
\sup_{z\in\D}(1-|z|^2)\left|\frac{zf''(z)}{f'(z)}\right|
\le 2a\Phi(L(f))
=\frac{2\delta}\pi \Phi(L(f)).
$$
Lemma \ref{lem:becker} now implies the required assertion.
\end{pf}

The above method also gives a norm estimate of the pre-Schwarzian derivative.
Though we do not use it in this note, we record it 
for the possible future reference.

\begin{prop}\label{prop:ps}
Suppose that $L(f)/l(f)<+\infty$ for a function $f\in\A.$
Then the pre-Schwarzian norm of $f$ is estimated as
$$
\|f\|\le \frac{2}\pi\left(1+L(f)\right)\log\frac{L(f)}{l(f)}.
$$
\end{prop}

\begin{pf}
Let $a=\frac1\pi\log \frac{L(f)}{l(f)}.$
By \eqref{eq:ps}, we have
$$
(1-|z|^2)\left|\frac{f''(z)}{f'(z)}\right|
\le 2a+2aL(f)(1-r^2)\frac{\arctanh r}r
$$
for $|z|=r<1.$
Since $(1-r^2)\arctanh r/r$ is decreasing in $0<r<1,$
the inequality $(1-r^2)\arctanh r/r\le 1$ holds.
Hence, we obtain
$$
(1-|z|^2)\left|\frac{f''(z)}{f'(z)}\right|\le 2a+2aL(f).
$$
\end{pf}

In order to prove Theorems \ref{thm:main1} and \ref{thm:main2},
the following technical result is helpful.
To state it, we introduce the auxiliary function
$$
H(x,c)=\frac{1-c}2 x+\frac{1+c}2 x\inv.
$$

\begin{lem}\label{lem:tech}
Let $c>1.$
If a number $x_1\in(0,1)$ satisfies the inequality 
$x_1\arctanh x_1<\frac{1+c}{2c},$
then $\Phi(c)<H(x_1,c).$
\end{lem}

\begin{pf}
Let $g(x)=x+c(1-x^2)\arctanh x.$
Then $g'(x)=1+c-2c x\,\arctanh x.$
Since $x\,\arctanh x$ (strictly) increases from $0$ to $+\infty$ when
$x$ moves from $0$ to $1,$ there exists a unique zero $x_0\in(0,1)$
of $g'(x)$ so that $g'(x)>0$ in $0<x<x_0$ and $g'(x)<0$ in $x_0<x<1.$
Note here that the assumption implies that $0<x_1<x_0.$
We see now that $g(x)$ takes its maximum at $x=x_0$ and therefore, we have
$$
\Phi(c)=g(x_0)=\frac{1-c}2 x_0+\frac{1+c}2 x_0\inv=H(x_0,c).
$$
Since $H_x(x,c)=(1-c)/2-(1+c)/(2x^2)<0,$ the function $H(x,c)$ is decreasing
in $x>0$ for a fixed $c>1.$
Hence, $x_1<x_0$ implies $H(x_0,c)<H(x_1,c),$ which proves the assertion.
\end{pf}

\begin{pf}[Proof of Theorem \ref{thm:main1}]
Let $\delta=\pi/3$ and set $c=e^{\delta/2}=e^{\pi/6}.$
If we take $x_1=17/22,$ then
$$
\frac{1+c}{2c}-x_1\arctanh x_1
=\frac{1+e^{-\pi/6}}2-\frac{17}{44}\log\frac{39}5
=0.00255\ldots>0.
$$
Therefore, Lemma \ref{lem:tech} yields
$$
\frac{2\delta}\pi\Phi(e^{\delta/2})
=\frac23\Phi(c)<\frac23 H(x_1,c)=\frac{773+195e^{\pi/6}}{1122}=0.982\ldots<1.
$$
We now apply Corollary \ref{cor:main} to obtain $\pi/3<\delta_0.$
\end{pf}

\begin{pf}[Proof of Theorem \ref{thm:main2}]
We will proceed in the same line as above.
Let $\delta=7\pi/25$ and set $c=e^\delta.$
We take $x_1=20/27$ and have
$$
\frac{1+c}{2c}-x_1\arctanh x_1=\frac{1+e^{-7\pi/25}}2-\frac{10}{27}\log\frac{47}7
=0.00219\ldots>0.
$$
Lemma \ref{lem:tech} now implies
$$
\frac{2\delta}\pi\Phi(e^{\delta})
=\frac{14}{25}\Phi(c)
<\frac{14}{25}H(x_1,c)=\frac{7903+2303e^{7\pi/25}}{13500}
=0.9965\ldots<1.
$$
We again apply Corollary \ref{cor:main} to obtain $7\pi/25<\delta_1.$
\end{pf}

\begin{rem}
We can slightly improve Theorems \ref{thm:main1} and \ref{thm:main2}
by changing the choice of $\delta$ and $x_1$ in the above proofs.
For instance, concerning Theorem \ref{thm:main1}, we can take 
$(\delta,x_1)=(\frac{22\pi}{65}, \frac{17}{22}), 
(\frac{87\pi}{257}, \frac{2765}{3578}),$
to have lower bounds $22\pi/65=1.06330\dots$ and
$87\pi/257=1.06349\dots,$ respectively, for $\delta_0.$
Numerical computations with Mathematica 8 suggest that
the solution to the equation $\frac{2\delta}\pi\Phi(e^{\delta/2})=1$
is about $\delta=1.0635213.$
Therefore, it seems that we would obtain at most this value as a lower bound
for $\delta_0$ by the above method.

Similarly, concerning Theorem \ref{thm:main2}, we can take
$(\delta,x_1)=(\frac{25\pi}{89}, \frac{622}{839}), 
(\frac{127\pi}{452}, \frac{321}{433}),$
to have lower bounds  $25\pi/89=0.882469\dots$ and
$127\pi/452=0.882704\dots,$ respectively, for $\delta_1.$

We see that the numerical solution to the equation
$\frac{2\delta}\pi\Phi(e^{\delta})=1$ is about $\delta=0.8827139.$
Therefore, the above method seems to give only a lower bound of $\delta_1$
not better than this value.
\end{rem}

\section{Obtaining upper bounds: non-univalence of a specific function}

We will provide upper bounds for $\delta_0$ by checking non-univalence
of the function $F_a\in\A$ defined by \eqref{eq:Fa} for a suitably chosen
positive constant $a.$
Since $F_a$ has no simple form to express, it is not easy to determine
its univalence.
In this note, we will observe its Grunsky coefficients to examine
univalence, whereas we used Gronwall's area theorem (or its refinement
by Prawitz) to see that $a\le 1$ is necessary for $F_a$ to be univalent.

Let $f\in\A.$
The Grunsky coefficients $c_{j,k}$ of $f$ are defined by the series
expansion
\begin{equation}\label{eq:expand}
\log\frac{f(z)-f(w)}{z-w}
=-\sum_{j,k=0}^\infty c_{j,k}z^jw^k
\end{equation}
in $|z|<\varepsilon, |w|<\varepsilon$ for a small enough $\varepsilon>0.$
We remark here that the obvious symmetry relation $c_{j,k}=c_{k,j}$ holds.
Note also that $c_{j,0}~(j=1,2,\dots)$ are the logarithmic coefficients
of $f(z)/z,$ in other words,
$-\log[f(z)/z]=c_{1,0}z+c_{2,0}z^2+\cdots$ as we can see by letting $w=0$ in
\eqref{eq:expand}.
Grunsky's theorem was strengthened by Pommerenke as in the following
(see \cite[Theorem 3.1]{Pom:univ}).

\begin{lem}\label{lem:gr}
Let $f\in\A$ and $\{c_{j,k}\}$ be its Grunsky coefficients.
If $f$ is univalent on $|z|<1$ then
$$
\sum_{m=1}^\infty m\left|\sum_{k=1}^n c_{m,k}t_k\right|^2
\le\sum_{m=1}^n\frac{|t_m|^2}m
$$
holds for arbitrary $n\ge1$ and $t_1,\dots,t_n\in\C.$
\end{lem}

We remark that the Grunsky coefficients are usually defined for
the function $g(\zeta)=1/f(1/\zeta).$
This change affects only the coefficients $c_{j,0}=c_{0,j},$
which do not involve the Grunsky inequalities.
See \cite{Hummel64} for more information.

From Lemma \ref{lem:gr}, the inequality
\begin{equation}\label{eq:G}
\sum_{m=1}^n m\left|\sum_{k=1}^n c_{m,k}t_k\right|^2
\le\sum_{m=1}^n\frac{|t_m|^2}m
\end{equation}
follows for every $n$ and $t_1,\dots,t_n\in\C.$
This implies that the Hermitian matrix
$G_f(n)=(\gamma_{j,k}^{(n)})$ of order $n$ is positive semi-definite;
in other words, $\mathbf{t}G_f(n)\mathbf{t}^*\ge0$ for any
$\mathbf{t}=(t_1, \dots, t_n)\in\C^n,$ where
$$
\gamma_{j,k}^{(n)}=\frac{\delta_{j,k}}j-\sum_{m=1}^nmc_{m,j}\overline{c_{m,k}}
\qquad (1\le j,k\le n),
$$
$\delta_{j,k}$ means Kronecker's delta and $\mathbf{t}^*$ is the conjugate
transpose of $\mathbf{t}$ as a matrix.

Letting $t_k=\delta_{j,k}$ in \eqref{eq:G},
we have $\sum_{m=1}^n m|c_{m,j}|^2\le 1/j$ for $j\le n,$  which implies
$|c_{m,j}|\le 1/\sqrt{mj}\le 1$ for $m,j\ge1.$
This guarantees that the series expansion in \eqref{eq:expand}
is convergent in $|z|<1, |w|<1,$ and therefore, that $f$ is univalent on $\D.$
We shall call $G_f(n)$ the {\it Grunsky matrix} of order $n.$
We have observed the following assertion.

\begin{cor}\label{cor:pom}
A function $f\in\A$ is univalent on $\D$ if and only if
its Grunsky matrix $G_f(n)$ of order $n$
is positive semi-definite for every $n\ge1.$
\end{cor}

In order to compute the Grunsky coefficients of $F_a(z),$
it is convenient to have recursion formulae for relavant coefficients.
The following elementary lemma gives a recursion formula for exponentiation.

\begin{lem}\label{lem:rec}
Let $g(z)=b_1z+b_2z^2+\cdots$ be a given function analytic around $z=0$
and let $h(z)=e^{g(z)}=c_0+c_1z+c_2z^2+\cdots.$
Then $c_n$ can be computed recursively by
$c_0=1$ and
$$
c_n=\frac1n\sum_{k=0}^{n-1}(n-k)b_{n-k}c_k\quad (n\ge1).
$$
\end{lem}

\begin{pf}
Compare the coefficients of the power series expansions
of both sides of $h'(z)=g'(z)h(z).$
\end{pf}

We turn to the function $F_a(z)$ for a fixed $a>0.$
In view of \eqref{eq:Q},
we see that the relation \eqref{eq:Fa} can also be expressed by
$zF_a'(z)/F_a(z)=Q_a(z)=\exp(2a\arctan z).$
In particular, the range of the function $zF_a'(z)/F_a(z)$ is
the annulus $e^{-\pi a/2}<|w|<e^{\pi a/2}$ and, in particular,
$l(F_a)=e^{-\pi a/2},~L(F_a)=e^{\pi a/2}$ and $L(F_a)/l(F_a)=e^{\pi a},$
as is already announced in Introduction.
Using the formula
$$
\arctan z=\sum_{n=0}^\infty\frac{(-1)^n}{2n+1}z^{2n+1}
$$
together with the last lemma,
we can compute the Taylor coefficients $b_n$ of $Q_a(z)$ recursively.
(See also \cite{Tod95} for additional information about the coefficients.)
In this way, we obtain
\begin{align*}
\frac{zF_a'(z)}{F_a(z)}&=Q_a(z)=\sum_{n=0}^\infty b_nz^n \\
&=1+2az+2a^2z^2+\frac23a(2a^2-1)z^3+\frac23a^2(a^2-2)z^4+\cdots.
\end{align*}
Dividing by $z$ and integrating it with respect to $z,$ we obtain
%$$
%\frac{F_a'(z)}{F_a(z)}=\frac1z+\sum_{n=1}^\infty b_nz^{n-1}
%$$
%leads to
$$
\log\frac{F_a(z)}{z}=\sum_{n=1}^\infty\frac{b_n}{n}z^n
=2az+a^2z^2+\frac29a(2a^2-1)z^3+\frac16a^2(a^2-2)z^4+\cdots.
$$
We again use Lemma \ref{lem:rec} to compute the Taylor coefficients
of $F_a(z)/z$ recursively and
finally arrive at the representation
\begin{align*}
F_a(z)&=z\exp\left(\sum_{n=1}^\infty\frac{b_n}{n}z^n\right) \\
&=z+2az^2+3a^2z^3+\frac29a(17a^2-1)z^4+\frac19a^2(38a^2-7)z^5+\cdots.
\end{align*}

In order to compute the Grunsky coefficients, we use the following formulae.
These formulae are essentially known.
See \cite{Jab66} and \cite[Formula (2.13)]{Curt71} for example.
However, since we could not find exactly the same formula in the literature,
we state it as a lemma with proof in this note.

\begin{lem}\label{lem:recg}
The Grunsky coefficients $c_{j,k}$ of a function $f(z)=z+a_2z^2+\cdots$
in $\A$ satisfy the recursion formulae
\begin{equation}\label{eq:grc}
c_{j,k}=\sum_{l=1}^{k-1}\frac lk a_{k-l}c_{j+1,l}-\sum_{m=1}^j a_{m+1}c_{j-m,k}
-\frac{a_{j+k+1}}k
\end{equation}
for $j\ge0$ and $k\ge1.$
\end{lem}

\begin{pf}
Differentiating both sides of \eqref{eq:expand} with respect to $w,$ 
we obtain the relation
$$
wf'(w)-w\frac{f(z)-f(w)}{z-w}
=(f(z)-f(w))\sum_{j,k=0}^\infty kc_{j,k}z^jw^k.
$$
Letting $a_1=1,$ we compute first the left-hand side of the relation:
$$
\text{(LHS)}=\sum_{n=1}^\infty a_n\big(nw^n-w(z^{n-1}+\cdots+zw^{n-2}+w^{n-1})\big)
=\sum_{n=1}^\infty a_n\big( (n-1)w^n-z^{n-1}w-\cdots-zw^{n-1}\big).
$$
On the other hand, 
\begin{align*}
\text{(RHS)}&=\sum_{n=1}^\infty\sum_{j,k=0}^\infty ka_nc_{j,k}(z^{j+n}w^k-z^jw^{k+n}) \\
&=\sum_{l,m=0}^\infty\left(
\sum_{n=1}^l ma_nc_{l-n,m}-\sum_{n=1}^m (m-n)a_nc_{l,m-n}
\right)z^lw^m.
\end{align*}
Comparing the coefficients of the term $z^lw^m,$ we get
$$%\begin{align*}
-a_{l+m}
%&=\sum_{n=1}^l ma_nc_{l-n,m}-\sum_{n=1}^m (m-n)a_nc_{l,m-n} \\
=mc_{l-1,m}+\sum_{n=2}^l ma_nc_{l-n,m}-\sum_{n=1}^m (m-n)a_nc_{l,m-n}
$$%\end{align*}
for $l\ge1$ and $m\ge 1.$
We now let $(j,k)=(l-1,m)$ to obtain the required relation.
\end{pf}

We can now compute $c_{j,k}$ recursively.
Indeed, first we apply \eqref{eq:grc} with $k=1$ to compute
$c_{j,1}$ recursively in $j\ge0:$
$$
c_{j,1}=-\sum_{m=1}^j a_{m+1}c_{j-m,1}-a_{j+2},\qquad j\ge0.
$$
If we determine $c_{l,m}$ for all $l\ge0$ and $1\le m<k,$ then
we use \eqref{eq:grc} to give $c_{j,k}$ recursively in $j\ge0.$
Practically, to determine $c_{j,k},$ it is enough to start with
$c_{l,1}$ for $0\le l\le j+k-1,$ which determine
$c_{l,2}$ for $0\le l\le j+k-2,$ and so on.
In this way, we can compute the Grunsky matrix $G(n)=G_{F_a}(n).$
For instance, $G(1)=[1-a^4]$ and
$$
G(2)=\frac1{81}\begin{bmatrix}
81-8a^2-97a^4-8a^6 & -14a^3(1+a^2)^2 \\
-14a^3(1+a^2)^2 & 81/2-4a^2-10a^4+10a^6-49a^8/2
\end{bmatrix}.
$$
We are now ready to give the upper bound in Theorems \ref{thm:main1}
and \ref{thm:main2}.

\begin{pf}[Computer-assisted proof of $\delta_0< 5\pi/7$]
We consider the Grunsky matrix $A_a=G(18)$ of order $18$ for the function
$f=F_{a}.$
We computed $A_a$ symbollically
with the help of Mathematica 8 but we will not
give a list of the elements of $A_a$ due to limitation of the space.
Let $a_0=5/7.$
We will show that $F_{a}$ is not univalent for $a$ close enough to $a_0.$

We see that $A_{a_0}$ is a square matrix of order 18 with rational elements.
Mathematica 8 can compute its eigenvalues and corresponding eigenvectors
numerically.
In this way, we found that one eigenvalue of $A_{a_0}$ was apparently negative.
Since numerical computations might not be reliable enough,
we will make this observation rigorous.
By approximating an eigenvector corresponding to the negative eigenvalue,
we find that the rational vector
$$
\mathbf{v}=\left(-\frac13, -\frac16, \frac3{10}, \frac3{10}, -\frac16,
-\frac13,0,\frac13,\frac16,-\frac15, -\frac15, \frac1{10},
\frac15, \frac1{10}, -\frac15, -\frac16, \frac15, \frac16\right)
$$
satisfies
$$
\mathbf{v}A_{a_0}\mathbf{v}^*
=-\frac{37\cdot 61\cdot 102353087\cdot 29977321169
\cdot N}
{3^{49}\cdot 5^{16}\cdot 7^{92}\cdot 11^{12}\cdot 13^4
\cdot 17^3\cdot 19^4 \cdot 23^4\cdot 29^2\cdot 31^4}<0.
$$
Here, $N=
76346348854682571404146112285557118341692971860401383400032365610149$
$904921555392748616477613599662190674795168801824208283713$ 
is an integer with 125 digits, which cannot be factorized anymore
by Mathematica 8.
Therefore, $A_{a_0}$ is not positive semi-definite.
Since $\mathbf{v}A_a \mathbf{v}^*<0$ still holds for 
$a$ close enough to $a_0,$
we have $a_*<a_0$ by Corollary \ref{cor:pom}, where $a_*$ is
the number defined in the Introduction.
We thus have seen that $\delta_1\le\delta_0\le \pi a_*<\pi a_0=5\pi/7.$
\end{pf}

\def\cprime{$'$} \def\cprime{$'$} \def\cprime{$'$}
\providecommand{\bysame}{\leavevmode\hbox to3em{\hrulefill}\thinspace}
\providecommand{\MR}{\relax\ifhmode\unskip\space\fi MR }
\providecommand{\MRhref}[2]{%
  \href{http://www.ams.org/mathscinet-getitem?mr=#1}{#2}
}
\providecommand{\href}[2]{#2}

%\bibliography{papers}
\end{document}